\documentclass[11pt, oneside]{article} 
\usepackage{latexsym}
\usepackage{amsfonts}
\usepackage{rawfonts}
\usepackage{graphicx}			
\usepackage{amsmath}										
\usepackage{amssymb}
\usepackage{geometry}                		
\input{prepictex}
\input{pictex}
\input{postpictex}

\newtheorem{thm}{Theorem}
\newtheorem{prop}[thm]{Proposition}
\newtheorem{cor}[thm]{Corollary}

\newtheorem{lem}[thm]{Lemma}

\newtheorem{pf}{Proof}
\newtheorem{rmk}{Remark}

\begin{document}

\title{On the Spectrum of weighted Laplacian operator and its application to uniqueness of K\"ahler Einstein metrics}

\author{Long Li}

\maketitle{} 
\begin{abstract}
The purpose of this paper is to provide a new proof of Bando-Mabuchi's uniqueness theorem of K\"ahler Einstein metrics on Fano manifolds, based on the convexity of $Ding$-functional on Chen's weak $\mathcal{C}^{1,\bar{1}}$ geodesic without
using any further regularities. Unlike the smooth case, the lack of regularities on the geodesic forbids us to use spectral formula of the weighed Laplacian operator directly. However, we can use smooth $\epsilon$-geodesics to approximate the weak one, then prove that a sequence of eigenfunctions will converge into the first eigenspace of the weighted Laplacian operator.
\end{abstract}
\section{Introduction}
The study of K\"ahler Einstein metrics on Fano manifolds is an old but lasting subject in complex geometry: on geometrical point of view, it characterizes the manifold with constant Ricci curvature, i.e. the K\"ahler metric satisfies 
\[
Ric( \omega) = \omega;
\]
on analytical point of view, the complex Monge-Amp\`ere equations arise from the study of this curvature equation, i.e. the K\"ahler potential $\varphi\in\mathcal{H}$ is the solution of the following equation
\[
(\omega_0 + i\partial\bar{\partial}\varphi)^{n} = e^{h - \varphi} \omega_0^n
\]
where $\mathcal{H}:= \{\omega = \omega_0 + i\partial\bar{\partial}\varphi> 0 \}$. Now as a PDE problem on manifolds, it's natural to ask two questions - existence and uniqueness. After Yau's celebrated work[14] on solving the Calabi Conjecture, Tian's $\alpha$ invariant[13] gives a sufficient condition to solve Monge Amp\`ere 
equation on Fano manifolds in 1980's. Then many people contribute to this problem during these years. And quite recently, Chen-Donalson-Sun's work([8], [9], [10]) proves the existence of K\"ahler Einstein metrics on Fano manfolds is equivalent to K-stability condition.This settles a long standing stability conjecture on K\"ahler Einstein metrics which goes back to Yau.
\\
\\
The problem of uniqueness of K\"ahler Einstein metrics on Fano manifolds also keeps attractive during these years. It is first proved by Bando and Mabuchi[1] in 1987, and we will give an alternative proof in this paper. The statement is as follows
\begin{thm}
Let $X$ be a compact complex manifold with $-K_X>0$. Suppose $\omega_1$ and $\omega_2$ are two K\"ahler Einstein metrics on $X$, then there is a holomorphic automorphism $F$, such that 
\[
F^*(\omega_2) = \omega_1
\] 
where this $F$ is generated by a holomorphic vector field $\mathcal{V}$ on X.
\end{thm}

They solve this problem by considering a special energy(Mabuchi energy) decreasing along certain continuity path. Then the existence of 
weak $\mathcal{C}^{1,\bar{1}}$ geodesic between any two smooth K\"ahler potentials is proved by X.X.Chen[7] in 2000, and this idea turns out to be an important tool in proving uniqueness theorems. For instance, Berman[3] gives a new proof of Bando-Mabuchi's theorem by arguing the geodesic connecting two K\"ahler Einstein metrics is actually smooth. And Berndtsson[5] proves the uniqueness of possible singular K\"ahler Einstein metrics along $\mathcal{C}^0$ geodesics. He observes the $Ding$-functional is convex along these geodesics from his curvature formula on the Bergman kernel[6]. Moreover, this curvature formula plays a major role to create a holomorphic vector fields when the functional is affine. This method is used by Berman again to prove the uniqueness of Donaldson's equation[2], and generalized to the $klt-pairs$ in [BBEGZ12]. 
\\
\\
The idea of this paper is also initiated from the convexity of $Ding$-functional along geodesics from a different perspective. However, instead of using Berndtsson's curvature formula, we are going to use the Futaki's formula(refer to Section 2) of weighted Laplacian operator to derive the holomorphic vector fields. Unlike the former case, here the main difficulty arises from the change of metrics during the convergence of Laplacian operators. Fortunately, we have control on the mixed derivatives $\partial_{\alpha}\partial_{\bar{\beta}}\phi$ on the product manifold, i.e. Chen's existence theorem of weak geodesic[7] guarantees a uniform bound of mixed second derivatives of the potential in both space and time directions on the geodesic. Moreover, we can perturb the weak geodesic to a sequence of nearby smooth metrics $\{g_{\epsilon}\}$ with mixed second derivatives under control[7].
\\
\\
However, this is not quite enough for our purpose, because we are lack of a uniform lower bound of these metrics, and the lower bound of metrics(or the upper bounds of the inverse metric) is crucially involved in the weighted Laplacian operator as 
\[
\Box_{\phi_g}u = \partial^{\phi_g} (\omega_g \lrcorner \bar{\partial} u)
\]
where $u $ is a smooth function on $X$. More fundamentally, it plays an important role in the $L^2$ norm of $(0,1)$ forms as 
\[
<\xi,\eta >_g = \int_X g^{i\bar{j}}\xi_{i}\overline{\eta_{j}} d\mu_g 
\]
where $\xi = \xi_{i}d\bar{z}^{i}$ and $\eta= \eta_{j}d\bar{z}^j$. This forbids us to use standard $L^2$ theorems to get the a uniform control. We will overcome this difficulty in Section 5 by separating the Laplacian equation to two equations, i.e.
\[
\omega_g\wedge v_g = \bar{\partial} u
\] 
and 
\[
\partial^{\phi_g} v = \Box_{\phi_g}u.
\]
This idea is initiated from solving the equation $\partial^{\phi}v= \pi_{\perp}\phi'$ in Berndtsson's work[5]. Then a crutial $W^{1,2}$ estimate of the sequence of vector fields $v_{g_{\epsilon}}$ shows it converges to some vector fields $v_{\infty}$ in strong $L^2$ norm, and a further $L^1$ estimate on $\bar{\partial}v_{g}$ indicates the vector fields $v_{\infty}$ is in fact holomorphic, under certain conditions(refer to proposition 12). This solves our problem on fiber direction, but on time direction we need to argue the holomorphic vector field keeps to be a constant. This is guaranteed since it corresponds to an eigenfuntion in the first eigenspace of the weighted Laplacian operator and satisfies the geodesic equation.
\\
\\
$\mathbf{Acknowledgement}$: I would like to express my great thanks to Prof. Xiuxiong Chen, who suggested me to do this problem and showed me the way of approach when the geodesic is smooth. I would also thank to Prof. Eric Bedford, Prof. Song Sun, Prof. Weiyong He, and Dr. Kai Zheng for helpful discussion. And especially, I want to thanks Prof. Futaki for pointing out one error in the old version. Finally, the suggestion from Chengjian Yao also helps me to make this paper more clear. 

\section{Futaki's formula and Hessian of $Ding$-functional}
The manifolds $X$ in our consideration is Fano, then we can assume the K\"ahler class $[\omega] = c_1(X)$, i.e. for each K\"ahler metric $\omega_g$, there exists a smooth function $F_g$ such that
\[
Ric(\omega_g) - \omega_g = i\partial\bar{\partial}F_g,
\]
hence we can define a weighted volume form as $e^F\det g$(we will write $F_g$ as $F$ when there is no confusion), and a pairing for any $u,v\in \mathcal{C}^{\infty}(X)$ 
\[
(u,v)_g = \int_X u \bar{v} e^{F}\det g,
\]
then Futaki[9] considers a weighted Laplacian operator 
\[
\Delta_{F} u = \Delta_g u - \nabla^j u \nabla_j F.
\]
the reason to do this is because the new Laplacian operator is easy to do integration by parts under the weighted volume form
\[
\int_X (\Delta_{F} u )\bar{u} e^{F} \det g = -\int_X (\nabla_j\nabla^j u + \nabla^j u \nabla_j F ) \bar{u} e^{F}\det g
\]
\[
= \int_X \nabla^j u \nabla_j\bar{u}e^{F}\det g
\]
\[
= \int_X |\bar{\partial} u|^2 e^{F} \det g
\]
where the norm of the 1-form is take with respect to the metric $g$. Hence it's an elliptic operator, and its spectral is discrete as $0 < \lambda_1 <\lambda_2 < \cdots$. Then for each eigenfunction $\Delta_{F} u = \lambda u$, Futaki[11] writes the following formula
\[
\lambda \int_X |\bar{\partial} u|^2 e^{F}\det g = \int_X |\bar{\partial} u |^2 e^{F}\det g + \int_X |L_g u|^2 e^{F}\det g
\]
where $L_g $ is a second order differential operator defined as
\[
L_g u = \nabla_{\bar{j}}\nabla^i u \frac{\partial}{\partial z^i}\otimes d\bar{z}^j.
\] 
Now observe the RHS of Futaki's formula is in fact $\int_X |\Delta_{F_g}u|^2 e^F\det g$, we can generalize it to all smooth function as
\begin{lem}
For any smooth function $u$ on $X$, we have
\[
\int_X | \Delta_{F}u|^2 e^F\det g =  \int_X |\bar{\partial} u |^2 e^{F}\det g + \int_X |L_g u|^2 e^{F}\det g.
\]
\end{lem}
\begin{pf}
we can decompose $u = \Sigma_0^{\infty} a_i(u) e_i $ into the eigenspace of the operator $\Delta_{F_g}$, and notice that the eigenfunction $e_i$ is orthogonal with respect to each other under the weighted volume form and metric $g$. Then the first two terms in above equation will preserve this orthogonality, i.e. choose eigenfunctions $u$ and $w$ of $\Delta_F$ which are orthogonal to each other, then 
\[
\int_X |\bar{\partial}u+ \bar{\partial}w |^2 e^F \det g = \int_X |\bar{\partial}u|^2 e^F \det g + \int_X |\bar{\partial}w|^2e^F\det g
\]
and 
\[
\int_X |\Delta_F u+ \Delta_F w |^2 e^F \det g =  \int_X |\Delta_F u|^2 e^F \det g + \int_X |\Delta_F w|^2e^F\det g
\]
Moreover, the differential operator $L_g$ keeps this orthogonality of eigenfunctions, but first notice
\[
F_{,\alpha\bar{\beta}} = R_{\alpha\bar{\beta}} - g_{\alpha\bar{\beta}}
\]
from the definition of $F$, then we compute as follows
\[
\int_X \langle L_g u, L_g w \rangle_g e^{F}\det g = \int_X g^{\alpha\bar{\lambda}}g^{\mu\bar{\beta}} u_{,\bar{\lambda}\bar{\beta}} \bar{w}_{,\mu\alpha}e^F\det g
\]
\[
= -\int_X g^{\alpha\bar{\lambda}}g^{\mu\bar{\beta}} u_{,\bar{\lambda}\bar{\beta}\alpha}\bar{w}_{,\mu} e^F\det g - \int_X g^{\alpha\bar{\lambda}}g^{\mu\bar{\beta}} u_{,\bar{\lambda}\bar{\beta}}\bar{w}_{,\mu}
F_{,\alpha} e^F\det g
\]
\[
= - \int_X g^{\alpha\bar{\lambda}}g^{\mu\bar{\beta}} u_{,\bar{\lambda}\alpha\bar{\beta}}\bar{w}_{,\mu}e^F\det g  
- \int_X g^{\mu\bar{\beta}} R_{\bar{\beta}}^{\bar{\gamma}} u_{,\bar{\gamma}}\bar{w}_{,\mu} e^F\det g 
\]
\[
+ \int_X g^{\alpha\bar{\lambda}}g^{\mu\bar{\beta}} u_{,\bar{\lambda}}\bar{ w}_{,\mu\bar{\beta}} F_{,\alpha} e^F\det g 
+ \int_X g^{\mu\bar{\beta}} u_{,\bar{\lambda}}\bar{w}_{,\mu} F^{, \bar{\lambda}}_{\bar{\beta}} e^F\det g 
+\int_X g^{\alpha\bar{\lambda}}g^{\mu\bar{\beta}} u_{,\bar{\lambda}}\bar{w}_{,\mu} F_{,\alpha} F_{,\bar{\beta}} e^F\det g
\]
\[
= \int_X g^{\alpha\bar{\lambda}}g^{\mu\bar{\beta}}u_{,\bar{\lambda}\alpha}\bar{w}_{,\mu\bar{\beta}} e^F\det g + \int_X g^{\alpha\bar{\lambda}}g^{\mu\bar{\beta}} u_{\bar{\lambda}\alpha}\bar{w}_{,\mu}
F_{,\bar{\beta}}e^F\det g
\]
\[
+ \int_X g^{\alpha\bar{\lambda}}g^{\mu\bar{\beta}} u_{,\bar{\lambda}}\bar{ w}_{,\mu\bar{\beta}} F_{,\alpha} e^F\det g 
+\int_X (g^{\alpha\bar{\lambda}} u_{,\bar{\lambda}}F_{,\alpha}  )(   g^{\mu\bar{\beta}} \bar{w}_{,\mu}  F_{,\bar{\beta}}) e^F\det g - \int_X g^{\mu\bar{\beta}}u_{,\bar{\beta}}\bar{w}_{\mu}e^F\det g
\]
\[
= \int_X ( g^{\alpha\bar{\lambda}}u_{,\alpha\bar{\lambda}} + g^{\alpha\bar{\lambda}}u_{,\bar{\lambda}}F_{,\alpha})(g^{\mu\bar{\beta}}\bar{w}_{\mu\bar{\beta}} + g^{\mu\bar{\beta}}\bar{w}_{\mu}F_{,\bar{\beta}}) e^F \det g
\]
\[
=\int_X (\Delta_F u, \Delta_F w )_g e^F \det g = 0.
\]
\end{pf}

Next let's consider an easy case: according to He[12], the second derivative of $Ding$-functional on a smooth geodesic equals
\[
\frac{\partial^2 \mathcal{D}}{\partial t^2} = (\int_X e^{F_g} \det g)^{-1} \{ \int_X ( |\bar{\partial}\varphi'|^2_g  - (\pi_{\perp}\varphi')^2)e^{F_g}\det g      \} 
\]
where the metric $g$ is induced by the K\"ahler form $\omega_{\varphi}$, and the projection operator is defined as $\pi_{\perp} u = u- \int_X u e^{F_g} \det g / \int_X e^{F_g}\det g $. This implies
$Ding$-functional is convex along smooth geodesics. Now suppose there is a smooth geodesic connecting two K\"ahler Einstein metrics, the $Ding$-functional must keep to be a constant along it. Hence we get
\[
\int_X  |\bar{\partial}\varphi'|^2_g e^{F_g}\det g  = \int_X (\pi_{\perp}\varphi')^2 e^{F_g}\det g,
\]
then we see the first eigenvalue $\lambda_1$ of the weighted Laplacian operator $\Delta_{F_g}$ is $1$, and $\pi_{\perp}\varphi'$ belong to the first eigenspace, i.e.
\[
\Delta_{F_g}(\pi_{\perp}\varphi') = \pi_{\perp}\varphi'.
\]
Now by Futaki's formula, we see 
\[
L_g (\pi_{\perp}\varphi') = 0,
\]
then the induced vector field $V _t = \nabla^i\varphi' \frac{\partial}{\partial z^i}$ is holomorphic on $X$. Moreover, let's differentiate this vector field with respect to $t$ on the geodesic
\[
(g^{j\bar{k}}\varphi'_{\bar{k}})' = g^{j\bar{k}}\varphi''_{\bar{k}} - g^{j\bar{q}} \varphi'_{p\bar{q}} g^{p\bar{k}}\varphi'_{\bar{k}}
\]
\[
= g^{j\bar{k}}(g^{\alpha\bar{\beta}}\varphi'_{\alpha}\varphi'_{\bar{\beta}})_{,\bar{k}} - g^{j\bar{q}} \varphi'_{p\bar{q}} g^{p\bar{k}}\varphi'_{\bar{k}}
\]
\[
= g^{j\bar{k}} g^{\alpha\bar{\beta}} \varphi'_{\alpha}\varphi'_{,\bar{\beta}\bar{k}} = 0
\]
by the holomorphicity of $V_t$. Finally, this gives us a holomorphic vector field $\mathcal{V} = V_t - \partial/\partial t$ on $X\times S$, and its induced automorphism will give the uniqueness of the two K\"ahler Einsteim metrics.

\section{Some $L^2$ theorems }

In this section, we are going to use $L^2$ theorem to investigate the weighted Laplacian operator $\Delta_{F_g}$ and its spectrum, then we shall project our target to the front eigenspace in the proof of uniqueness theorem. First notice that we always have $\lambda_1 \geqslant 1$ by Futaki's formula. Then we are going to introduce some notations.
\\
\\
From now on, we shall assume the manifold $X$ admits non-trivial holomorphic vector fields, and $H^{0,1}(X) =0$. Then fix one $t$ and restrict our attention to this fiber $X\times \{t \}$. Since $-K_X = [\omega]$, we can write
\[
\omega_g = i\partial\bar{\partial}\phi_g
\]
where $\phi_g$ is a plurisubharmonic metric on the line bundle $-K_X$. We claim the measure
\[
e^{F_g}\det g = e^{-\phi_g},
\]
and this is because locally $F_g = -\log\det g - \phi_g$. Then naturally the pairing between functions on $X$ with this weight can be written as
\[
(u,v)_g = \int_X u\bar{v} e^{-\phi_g}.
\]
Here is the $L^2$ theorem coming to play with. Let's consider the space of all $L^2$ bounded $-K_X$ valued $(n,0)$ forms under the metric $\phi_g$, i.e. it consists of every function $u$ on $X$ such that
\[
 \int_X  |u |^2 e^{-\phi_g} < +\infty,
\]
we denote this space as $L^2_{(n,0)}(-K_X, \phi_g)$, and similarly we can consider all $L^2$ bounded $-K_X$ valued $(n, 1)$ forms under the weighted norm
\[
\int_X g^{\alpha\bar{\beta}}v_{\alpha}\overline{v_{\beta}} e^{-\phi_g} < +\infty,
\] 
and we denote this space as $L^2_{(n,1)}(-K_X, \phi_g)$, then we can define an unbounded operator $\bar{\partial}$ between them 
\[
\bar{\partial}: L^2_{(n,0)}(-K_X, \phi_g) \dashrightarrow L^2_{(n,1)}(-K_X, \phi_g).
\]
Notice that the domains of these two operator are not the whole $L^2$ spaces. In fact, we can define 
\[
dom(\bar{\partial}): = \{ u\in  L^2_{(n,0)}(-K_X, \phi_g) ; \ \bar{\partial}u\in     L^2_{(n,1)}(-K_X, \phi_g)  \},
\]
but it is not densely defined in $L^2$ space when $g_{\phi}$ is a  $\mathcal{C}^{1,\bar{1}}$ solution of geodesic equation on a fiber $X\times\{ t\}$. Hence we should consider the Hilbert space $\mathcal{H}_1$ to be the closure of $dom(\bar{\partial})$ in $L^2_{(n,0)}(-K_X, \phi_g)$. We claim that $\mathcal{H}_1$ is not empty.\\
\\
First notice that for any non-trivial holomorphic vector field $v \in L^2_{(n-1,0)}(-K_X)$, we can solve the following equation
\[
\bar{\partial}u = \omega_{g_{\phi}}\wedge v,
\]
since $\bar{\partial}( \omega_{g_{\phi}}\wedge v) = 0$ in the sense of distributions, but $\ker(\bar{\partial}) = Range(\bar{\partial})$ from $H^{0,1}(X) =0$. Next, consider the subspace 
$\mathcal{W}$ containing all such $u$, i.e. define
\[
\mathcal{W}: = \{u\in L^2_{(n,0)}(-K_X, \phi_g);\ \bar{\partial}u =\omega_{g_{\phi}}\wedge v, \   \forall v\in L^2_{(n-1,0)}(-K_X) \},
\]
then it is a non-empty subspace in $L^2_{(n,0)}(-K_X, \phi_g)$, and it's easy to check 
\[
\mathcal{W} \subset dom(\bar{\partial}),
\]
hence we proved the claim. Now $\bar{\partial}$ is a densely defined, closed operator on the Hilbert space $\mathcal{H}_1$ - it's closed from the continuity property of differential operators in the distribution sense. We can discuss its Hilbert adjoint operator $\bar{\partial}^*_{\phi}$, which is a densely defined, closed operator on $L^2_{(n,1)}(-K_X, \phi_g)$. Moreover, they have closed ranges 
\begin{lem}
$\bar{\partial}$ and $\bar{\partial}^{*}_{\phi_g}$ are densely defined, closed operators with closed ranges.
\end{lem}
\begin{pf}
We need to estimate the $L^2$ norm of $\bar{\partial}u$. Take $h$ to be a fixed smooth metric with positive Ricci curvature on $X$, and $u\in dom(\bar{\partial})\cap \ker (\bar{\partial})^{\perp}$, we have
\[
\int_X |\bar{\partial}u|_g^2 e^{-\phi_g} \geqslant \int_X |\bar{\partial}u|_h^2 \det h 
\]
\[
\geqslant c \int_X |u|^2 \det h 
\]
\[
\geqslant c' \int_X |u|_g^2 e^{-\phi_g}.
\]
this estimate implies $\bar{\partial}$ has closed range, and hence its adjoint $\bar{\partial}^*_{\phi_g}$ by functional analysis reason.
\end{pf}
Then we can define the Laplacian operator as $\Box_{\phi_g} = \bar{\partial}^*_{\phi_g}\bar{\partial}$, where also as an unbounded closed operator, i.e.
\[
\Box_{\phi_g}: L^2_{(n,0)}(-K_X, \phi_g)\dashrightarrow L^2_{(n,0)}(-K_X, \phi_g)
\]
and its domain of definition is
\[
dom(\Box_{\phi_g}): =\{ u\in  L^2_{(n,0)}(-K_X, \phi_g);\  u\in dom(\bar{\partial})\  and \ \bar{\partial}u \in dom(\bar{\partial}^*_{\phi_g})   \}.
\]
we claim this operator also has closed range. and 
\begin{prop}
we have 
\[
\ker \Box_{\phi_g} = coker \Box_{\phi_g},
\]
hence they are both finite dimensional. 
\end{prop}
\begin{pf}
First note $\ker \Box_{\phi_g} = \ker \bar{\partial}$ is the 1 dimensional space of constant functions on $X$. In order to prove $coker \Box_{\phi_g}$ also has finite rank, it's enough to prove the weighted Laplacian operator has closed range, since it's self-adjoint
\[
coker \Box_{\phi_g} = R(\Box_{\phi_g})^{\perp} = \ker\Box_{\phi_g}.
\]
Now we are going to prove the closed range property, but this follows from the following estimate for $u\in dom(\Box_{\phi_g})\cap  \ker (\bar{\partial})^{\perp}$
\[
|| u ||_g^2 \leqslant  C ||\bar{\partial}u ||_g^2 
\]
\[
\leqslant C (\Box_{\phi_g} u, u)_g
\]
\[
\leqslant  2C || \Box_{\phi_g} u ||^2_g + \frac{1}{2} || u ||^2_g
\]
and hence 
\[
|| u ||^2_g \leqslant C' ||\Box_{\phi_g} u ||^2_g,
\]
which implies the claim.

\end{pf}
Notice that this is not enough to guarantee the existence of discrete spectral, but we have a further estimate,
\begin{lem}
For all $u\in dom(\Box_{\phi_g})\cap \ker(\bar{\partial})^{\perp}$, there is an uniform constant $C$, such that
\[
|| u ||_{W^{1,2}} \leqslant C || \Box_{\phi_g} u ||^2_g.
\]
\end{lem}
\begin{pf}
we still compare it with some fixed smooth weight(metric) $h$,
\[
||\bar{\partial}u ||_h^2 \leqslant C ||\bar{\partial}u ||_g^2
\]
\[
=C (\Box_{\phi_g}u, u )_g
\]
\[
\leqslant C ||\Box_{\phi_g} u||_g ||u||_g
\]
\[
\leqslant C' ||\Box_{\phi_g} u||_g ||u||_h
\]
\[
\leqslant C'' ||\Box_{\phi_g}u||_g ||\bar{\partial}u ||_h,
\]
then 
\[
||\bar{\partial}u ||_h^2  \leqslant C'' ||\Box_{\phi_g}u||_g.
\]
finally, an integration by part gives the desired estimate since
\[
\int_X h^{\alpha\bar{\beta}} u_{,\alpha}\overline{u_{,\beta}} \det h = -\int_X h^{\alpha\bar{\beta}} u_{,\alpha\bar{\beta}} \bar{u}\det h
\]
\[
= -\int_X h^{\alpha\bar{\beta}} u_{,\bar{\beta}\alpha} \bar{u}\det h
\]
\[
= \int_X h^{\alpha\bar{\beta}}u_{,\bar{\beta}}\overline{u_{,\bar{\alpha}}} \det h
\]
\end{pf}
Then we can discuss the spectral of $\Box_{\phi_g}$, when $g_{\phi}$ is the $\mathcal{C}^{1,\bar{1}}$ function. Suppose $\lambda$ is an eigenvalue of $\Box_{\phi_g}$, and 
let $\Lambda$ be the corresponding eigenspace, we claim
\begin{prop}
$\dim \Lambda < +\infty$
\end{prop}
\begin{pf}
Let $v_i\in \Lambda$ be a sequence of eigenfunctions with bound $L^2$ norm, i.e. $||v_i ||^2_g =1$, then since 
\[
|| v_i ||_{W^{1,2}} \leqslant C || \Box_{\phi_g} v_i ||_g
\]
\[
= C \lambda,
\]
hence there exists a $W^{1,2}$ function $v_{\infty}$ such that $v_i\rightarrow v_{\infty}$ in strong $L^2$ norm, by compact embedding theorem. And since 
$\Lambda =\ker (\Box_{\phi_g}- \lambda I)$ is a closed subspace of $L^2$
\[
v_{\infty}\in \Lambda.
\]
This implies every bounded sequence in $\Lambda$ has a convergent subsequence, i.e. the unit ball in $\Lambda$ is compact, hence $\dim\Lambda$ is finite. 
\end{pf}

Next we are going to discuss some computations when the weight $\phi_g$ is at least $\mathcal{C}^2$. First notice that formally
\[
< \Box_{\phi_g} u, v >_g \   = \ < \bar{ \partial} u, \bar{\partial} v>_g
\]
for any pairing $u, v$. It's easy to see
\[
\Box_{\phi_g}u = \Delta_{\phi_g}u
\]
for all smooth functions $u$, when the metric $\phi_g$ is smooth. If we look closer at these operators, there is a more computable way to express them. For this purpose, let's assume $\phi_g$ is a $\mathcal{C}^2$ metric, then for any $(n,1)$ form $\alpha $ with value in $-K_X$, 
\[
\bar{\partial}^*_{\phi_g} \alpha = \partial^{\phi_g} (\omega_g \lrcorner \alpha)
\]
where $\partial^{\phi} v = e^{\phi}\partial (e^{-\phi} v ) = \partial v - \partial\phi\wedge v$ for any $(n-1,0)$ form with value in $-K_X$(that is a vector field on $X$). Hence if we define 
\[
v= \omega_g\lrcorner \alpha,
\]
we will have 
\[
\bar{\partial}^*_{\phi_g} \alpha = \partial^{\phi_g} v
\]
and the weighted Laplacian operator could be computed as 
\[
\Box_{\phi_g} u = \partial^{\phi_g} (\omega_g\lrcorner \bar{\partial} u)
\]
for $u\in dom \Box_{\phi_g}\cap L^2_{(n,0)}(-K_X, \phi_g)$. Notice that there is commutation relation between the new defined operator $\partial^{\phi}$ and $\bar{\partial}$, that is
\begin{equation}
\partial^{\phi}\bar{\partial} + \bar{\partial} \partial^{\phi} = i\partial\bar{\partial}\phi\wedge\cdot
\end{equation}
Now if $u$ is any eigenfunction of the weighted Laplacian operator with eigenvalue $\lambda$, i.e. $\Box_{\phi_g} u = \lambda u$, we can decompose it into two equations 
\[
\omega_g \lrcorner \bar{\partial}u   = v \ \ \  \partial^{\phi_g} v = \lambda u.
\]
here we can write $v =  X\lrcorner 1$, where the constant function $1$ is read as an $(n,0)$ form with value in $-K_X$, and $X= X^{\alpha}\frac{\partial}{\partial z^{\alpha}}$ is a vector field in 
$(1,0)$ direction on the manifolds. Next we are going to prove Futaki's formula by the commutation equality.
\begin{lem}(Futaki's formula)
Let $u$ be a eigenfunction of weighted Laplacian with eigenvalue $\lambda$, i.e. $\Box_{\phi_g}u = \lambda u$, then 
\[
\lambda \int_X |\bar{\partial}u|^2_g e^{-\phi_g} = \int_X (|L_g u|^2 + |\bar{\partial}u|^2_g)e^{-\phi_g}.
\]
\end{lem}
\begin{pf}
First notice $u$ is pure real or imaginary. Hence here we will give the proof when $u$ is real valued - the case when $u$ is pure imaginary is similar. Now  by the commutation relation of $\partial^{\phi_g}$, we compute $\bar{\partial}(\lambda u)$
\[
-\partial^{\phi_g}\bar{\partial} v + i\partial\bar{\partial}\phi_g \wedge v = \lambda\bar{ \partial} u,
\]
notice that $i\partial\bar{\partial} \phi_g = \omega_g$, hence 
\[
-\partial^{\phi_g}\bar{\partial} v = (\lambda -1)\bar{\partial} u,
\]
pair it with $\bar{\partial}u$, 
\[
(\lambda -1) \int_X |\bar{\partial}u |^2_ge^{-\phi_g}=-\int_X \langle \partial^{\phi_g}\bar{\partial} v, \bar{\partial}u \rangle_g e^{-\phi_g} 
\]
\[
= \int_X -g^{\lambda\bar{\mu}}\partial_{\alpha} (e^{-\phi_g} \partial_{\bar{\mu}} X^{\alpha}) \overline {\partial_{\bar{\lambda}}u}
\]
\[
=\int_X \partial_{\bar{\mu}}X^{\alpha} \overline{\partial_{\bar{\alpha}}X^{\mu}} e^{-\phi_g}.
\]
Now notice that $X^{\alpha} = g^{\alpha\bar{\beta}}u_{,\bar{\beta}}$, under the normal coordinate when $g_{i\bar{j}} = \delta_{ij} \Lambda_i$,
\[
\partial_{\bar{\mu}}X^{\alpha}\partial_{\alpha}X^{\bar{\mu}} = g^{\alpha\bar{\beta} }u_{,\bar{\beta}\bar{\mu}} g^{\lambda\bar{\mu}}u_{, \lambda\alpha}
\]
\[
= \frac{1}{\Lambda_{\alpha} \Lambda_{\lambda}}u_{,\bar{\alpha}\bar{\lambda}} u_{,\lambda\alpha}
\]
\[
= \frac{1}{\Lambda_{\alpha}\Lambda_{\lambda}} u_{,\bar{\alpha}\bar{\lambda}} u_{,\alpha\lambda}
\]
\[
= g_{\alpha\bar{\beta}}g^{\lambda\bar{\mu}}  \partial_{\bar{\mu}}X^{\alpha} \overline{\partial_{\bar{\lambda}}X^{\beta}},
\]
hence we proved the Futaki's formula
\[
(\lambda -1)\int_{X}|\bar{\partial}u|^2_ge^{-\phi_g} = \int_X |\bar{\partial}X|^2_g e^{-\phi_g}.
\]
\end{pf}

\section{$Ding$-functionals along the approximation geodesics}
Let $X$ be an $n$ dimensional compact complex K\"ahler manifold with K\"ahler metric $\omega$, then we can write the K\"ahler form locally as 
\[
\omega = g_{\alpha\bar{\beta}}dz^{\alpha}\wedge d\bar{z}^{\beta}
\]
where $\alpha,\beta = 1,\cdots,n$. Take $S$ to be a cylinder, and $z^{n+1} = t + \sqrt{-1}s$ be its coordinate. Then $z = (z^1,\cdots, z^n, z^{n+1})$ is a point in $X\times S$, and we can define 
\[
\tilde{\omega} = g_{\alpha\bar{\beta}}dz^{\alpha}\wedge d\bar{z}^{\beta} + dz^{n+1}\wedge d\bar{z}^{n+1}
\]
as a K\"ahler metric on $X\times S$. And $\tilde{\varphi} = \varphi - |z^{n+1}|^2$ as the new potential. We shall write $\tilde{\omega}$ as $\omega$ and $\tilde{\varphi}$ as $\varphi$ when there is no fusion. Then Chen[7] proves the following two theorems

\begin{thm}(Existence of weak geodesic)
Let $\varphi_0, \varphi_1\in \mathcal{H}$, then there exists a unique $C^{1,\bar{1}}$ geodesic connecting them, i.e. the following homogenous Monge-Amp\`ere equation has a unique weak solution $\varphi\in \overline{\mathcal{H}}$(the closure is taken under the $C^{1,\bar{1}}$ topology) on $X\times S$ 
\[
\det (g_{i\bar{j}}+ \partial_{i}\partial_{\bar{j}}\varphi)_{(n+1)(n+1)} = 0
\]
where $i,j =1,\cdots, n+1$, and on the boundary $\partial(X\times S)$
\[
\varphi(0,s,z) = \varphi_0(z),\ \varphi(1,s,z) =\varphi_1(z)
\]
with the following estimate
\[
||\varphi||_{\mathcal{C}^1(X\times S)}+ \max\{|\partial_i\partial_{\bar{j}}\varphi    |  \} < C
\]
where $C$ is a uniform constant only depending on $\varphi_0$ and $\varphi_1$.
\end{thm}

\begin{thm}($\epsilon$- approximation geodeiscs)
Given $\varphi_0,\varphi_1\in\mathcal{H}$, we can have a sequence of approximation geodesics $\varphi_{\epsilon}(t, z)$ as follows: for each small $\epsilon>0$, there exists a unique solution of the equation 
\[
(\varphi_{tt} - |\partial_X \varphi'|_{g_{\varphi}}^2)\det(g_{\varphi}) = \epsilon \det h
\]
such that there exists a uniform constant $C$ with
\[
|\varphi'_t| + |\varphi''_t| +   |\varphi|_{\mathcal{C}^1} + \max\{  |\partial_{\alpha}\partial_{\bar{\beta}}\varphi| \} < C,
\]
and $\varphi_{\epsilon}$ converges to the $\mathcal{C}^{1,\bar{1}}$ geodesic $\varphi$ in the weak $\mathcal{C}^{1,\bar{1}}$ topology.
\end{thm}

Notice that for any plurisubharmonic metric $\phi$ on $-K_X$, we can write its potential as $\varphi = \phi - \phi_0$, where $\phi$ and $\phi_0$ are corresponding metrics on the line bundle 
$-K_X$. Now suppose $\phi_0, \phi_1$ are two smooth K\"ahler Einstein metrics on $X$, with their K\"ahler forms $\omega_i=i\partial\bar{\partial}\phi_i, i =0,1$ satisfying
\[
\omega_i^n = \frac{e^{-\phi_i}}{\int_X e^{-\phi_i}}.
\]
define the following functionals
\[
\mathcal{F}(\phi):= -\log\int_X e^{-\phi}
\]
and 
\[
\mathcal{E}(\phi):= \frac{1}{n+1} \Sigma_{j=0}^n \int_X \varphi \omega_0^j\wedge\omega_{\phi}^{n-j}
\]
where $\omega_{\phi} = i\partial\bar{\partial}\phi$. Then the $Ding$-functional is defined as 
\[
\mathcal{D} = -\mathcal{E} + \mathcal{F} = - \frac{1}{n+1} \Sigma_{j=0}^n \int_X (\phi-\phi_0) \omega_0^j\wedge\omega_{\phi}^{n-j}-\log\int_X e^{-\phi}.
\] 
Notice the along a curve of metrics $\phi_t$, the derivative of $Ding$-functional is
\[
\frac{\partial\mathcal{D}}{\partial t} =  \int_X \phi' ( - \omega_{\phi}^n + \frac{ e^{-\phi}}{\int_X e^{-\phi}}).
\]
we see the critical point of this functional is the K\"ahler Einstein metric, and its second derivative is 
\[
\frac{\partial^2 \mathcal{D}}{\partial t^2} = - \int_X (\phi'' - |\partial \phi'|_g^2 )\omega_{\phi}^n+ (\int_X e^{-\phi})^{-1} \{ \int_X(\phi'' - |\partial\phi'|^2_g )e^{-\phi}
+\int_X (|\partial\phi'|_g^2 - (\pi_{\perp}\phi')^2)e^{-\phi}   \}
\]
where the metric $g =i\partial\bar{\partial}\phi_t$, and if we denote the term $f = \phi'' - |\partial \phi'|_g^2 $, $c_t = \int_X e^{-\phi}$ and $\delta_t = |\partial\phi'|_g^2 - (\pi_{\perp}\phi')^2$, the equation reads 
\[
\frac{\partial^2 \mathcal{D}}{\partial t^2} = -\int_X f \omega_{\phi}^n + \int_X ( f+\delta_t) e^{-\phi} /c_t,
\]
then we are going to consider the behavior of $Ding$-functional on the approximation geodesic. First from Chen's theorem, we can find a $\mathcal{C}^{1,\bar{1}}$ geodesic $\phi_t$ connecting the two K\"ahler Einstein metrics. Moreover for any small $\epsilon>0$, there is the smooth approximation geodesic $\phi_{\epsilon}(t,z)$ connecting the two end points $\phi_0, \phi_1$, which converges weakly to the $\mathcal{C}^{1,\bar{1}}$ geodesic. Now if we consider the $Ding$-functional on these approximation geodesics, we have estimates
\[
\frac{\partial^2 \mathcal{D}}{\partial t^2} \geqslant -\epsilon \int_X \det h
\] 
from $f = \epsilon \det h / \det g > 0$ and $\int_X \delta_t e^{-\phi} > 0 $. Let $\epsilon\rightarrow 0$, we see that $Ding$-functional keeps to be convex on $\mathcal{C}^{1,\bar{1}}$ geodesic.
Now we can integrate it back along $t$ 
\[
\frac{\partial \mathcal{D}}{\partial t}(1) -\frac{\partial \mathcal{D}}{\partial t}(0) = \int_{X\times I} -f\omega^n_{\phi}dt + \int_{X\times I} f e^{-\phi}/c_t dt + \int_{X\times I} \delta_t e^{-\phi}/c_t dt,
\]
notice that at end points $\phi_0,\phi_1$ are both K\"ahler Einstein, hence the first derivative of $Ding$-functionals vanish. And on the approximation geodesic, we have the equation
\[
f \det g = \epsilon \det h
\] 
and $f \leqslant \phi'' < C$ uniformly independent of $\epsilon$. Then the equation above reads
\[
A \epsilon =    \int_{X\times I} f e^{-\phi}/c_t dt + \int_{X\times I} \delta_t e^{-\phi}/c_t dt 
\]
\[
\geqslant  \int_{X\times I} f e^{-\phi}dt + \int_{X\times I} \delta_t e^{-\phi}dt,
\]
because we have uniform $\mathcal{C}^0$ estimate on $\phi_{\epsilon}$. Now since we want to discuss the eigenfunctions on each fiber, we need to a lemma to pull back the estimate to fibers.
\begin{lem}
Suppose $F_{\epsilon}(t)$ is a sequence of non-negative function on $[0,1]$, with integration estimate
\[
\int_0^1 F_{\epsilon} dt < A \epsilon,
\]
then for almost everywhere $t\in[0,1]$, we can find a subsequence(depending on $t$) $F_{\epsilon_j}$, such that 
\[
F_{\epsilon_j} < C_t\epsilon_j
\]
where $C_t$ is a constant independent of $\epsilon$.
\end{lem}
\begin{pf}
Let $\tilde{F}_{\epsilon} = F_{\epsilon}/ \epsilon$, then by Fatou's lemma
\[
\int_0^1 \liminf_{\epsilon} \tilde{F}_{\epsilon} dt \leqslant \liminf_{\epsilon} \int_0^1 \tilde{F}_{\epsilon}dt \leqslant A,
\]
hence the function $\liminf_{\epsilon}\tilde{F}_{\epsilon} \in L^1$, i.e. for almost everywhere $t$, there is a subsequence $\tilde{F}_{\epsilon_j}$ and a constant $C_t$ such that 
\[
\tilde{F}_{\epsilon_j} < C_t,
\]
hence
\[
F_{\epsilon_j} < C_t \epsilon_j.
\]
\end{pf}
Now put $F_{\epsilon} = \int_X f_{\epsilon}e^{-\phi_{\epsilon}} + \int_X\delta_{\epsilon}e^{-\phi_{\epsilon}}$ and notice the two terms on RHS are both non-negative, we have proved
\begin{prop}
Consider the approximation geodesic $\phi_{\epsilon} $ connection two K\"ahler Einstein metrics. For almost everywhere $t$, there is a constant $C_t$, such that for each such $t$, there
exists a subsequence $\epsilon_j$, such that the following estimates
\[
\int_X f e^{-\phi} (\epsilon_j) < C_t \epsilon_j
\]
and 
\[
\int_X (|\partial\phi'|^2_g - (\pi_{\perp}\phi')^2)e^{-\phi} (\epsilon_j) < C_t \epsilon_j
\]
hold simutaneouly.
\end{prop}

\section{Convergence in the first eigenspace}
In this section, we shall focus our attention to the one fiber $X\times\{t \}$, and picked up a subsequence $\phi_{\epsilon_j}$ from above section. Then we can consider the sequence of weighted Laplacian operator $\Box_{\phi_{\epsilon}}$(we shall omit the subindex $j$ here). For each $\epsilon$, we can arrange its eigenvalues as $0<\lambda_1^{\epsilon}\leqslant \lambda_2^{\epsilon}\leqslant\cdots$, corresponding with one eigenfunction $e_i(\epsilon)$, i.e.
\[
\Box_{\phi_{\epsilon}} e_i(\epsilon) = \lambda_i^{\epsilon} e_i(\epsilon).
\]
Then let $u_{\epsilon}(z)$ be a sequence of smooth functions on $X$, such that $u_{\epsilon}\perp \ker\bar{\partial}$. Then it decomposes into the eigenspace of weighted Laplacian operator $\Box_{\phi_{\epsilon}}$, i.e.
\[
u_{\epsilon} = \Sigma_{i=1}^{N_{\epsilon}} a_i(\epsilon) e_{i}(\epsilon)
\] 
where $e_i\in \Lambda_i$, and in prior, $N_{\epsilon}$ could equal to $+\infty$ in the above notation.
Then we can consider the action by the weighted Laplacian operator on this sequence of functions, i.e. we can write $\Box_{\phi_{\epsilon}} u_{\epsilon}$ as
\[
v_{\epsilon} = \omega_{g_{\epsilon}}\lrcorner \bar{\partial}u
\]
and
\[
\partial^{\phi_{\epsilon}} v_{\epsilon} = \Sigma_{i=1}^{N_{\epsilon}} \lambda^{\epsilon}_i a_i(\epsilon)e_i(\epsilon).
\]
Under certain constraint, we claim these vector fields $v_{\epsilon}$ will converge to a holomorphic one with the same equation satisfied,
\begin{prop}
Let $u_{\epsilon}$ be a sequence of functions as above. Suppose it satisfies the following conditions:

1) $\Sigma_{i=1}^{N_{\epsilon}} |a_i(\epsilon)|^2 <  A$ for an uniform constant $A$, and the sums does not converge to zero. 

2) there exists a uniform constant $K$, such that $\lambda_{N_{\epsilon}}^{\epsilon} < K$ for each $\epsilon$

3) the following estimate holds
\begin{equation}
\int_X (|\bar{\partial} u_{\epsilon}|^2_{g_{\epsilon}} - (\pi_{\perp}u_{\epsilon})^2)e^{-\phi_{\epsilon}} < C \epsilon.
\end{equation}

then by passing to a subsequence, we have 
\[
u_{\epsilon}\rightarrow u_{\infty}
\]
in strong $L^2$ sense, where $u_{\infty}\in W^{1,2}$ is nontrivial. Moreover there exists a nontrivial holomorphic
$(n-1,0)$ form $v_{\infty}$ with value in $-K_X$, such that
\[
v_{\epsilon} \rightarrow v_{\infty}
\]
in strong $L^2$ sense, and the equation
\[
\omega_g \wedge v_{\infty} = \bar{\partial} u_{\infty}
\] 
holds in the sense of $L^2$ functions, where $g$ is the metric found on the $\mathcal{C}^{1,\bar{1}}$ geodesic.
\end{prop}
before proving the proposition, we need a lemma
\begin{lem}
Let $f_j, g_j$ be two sequence of $L^2$ functions with $|| f_j g_j   ||_{L^p} < C$ for some $p \geqslant 1$. Suppose that  $\int_X | f_j |^2 d\mu< C'$ and $g_j\rightarrow g\in L^2$ in $L^2$ norm, then there exists an $L^2$ function $f$ such that
\[
f_j g_j\rightarrow fg \in L^p
\]
in the sense of distributions.
\end{lem}
\begin{pf}
First note there exists an $L^2$ function $f$ such that $f_j\rightarrow f$ in weak $L^2$ topology. Then we check
\[
 \int_X ( fg -  f_j g_j ) d\mu =  \int_X g(f-f_j)d\mu + \int_X f_j(g-g_j)d\mu,
 \]
the first term on the RHS of above equation converges to zero from the weak convergence of $f_j$, and the second term converges to zero too, since 
\[
|\int_X f (g-g_j)d\mu|^2 \leqslant (\int_X |f|^2d\mu )(\int_X |g-g_j|^2d\mu) \rightarrow 0.
\]
hence $f_jg_j$ converges to $fg$ in the sense of distributions. Moreover, from the $L^p$ bound of $f_jg_j$, we have an $L^p$ function $k$ such that $f_jg_j\rightarrow k$ in weak $L^p$ topology. Then 
\[
fg = k
\]
as $L^p$ functions.
\end{pf}
\begin{rmk}
Suppose the sequence $|f_j|$ is uniformly bounded in lemma 13, then the limit $f$ is an $L^{\infty}$ function, then $fg\in L^2$ automatically. 
\end{rmk}
\begin{pf}(of proposition 12)
First we can write equation (2) as 
\[
\Sigma_{i=1}^{N_{\epsilon}} (\lambda^{\epsilon}_i  -1 )|a_i(\epsilon)|^2 < C\epsilon
\]
by Futaki's formula, we know 
\[
\int_X |L_{g_{\epsilon}} u_{\epsilon}|^2 e^{-\phi_{\epsilon}} = \Sigma_{i=1}^{N_{\epsilon}} \lambda_{i}^{\epsilon}(\lambda_i^{\epsilon} -1)|a_i(\epsilon)|^2
\]
\[
\leqslant KC\epsilon
\]
from condition (2) and (3). But if we write $v_{\epsilon} = X_{\epsilon}\lrcorner 1$ for some vector field $X_{\epsilon} = X_{\epsilon}^{\alpha}\frac{\partial }{\partial z^{\alpha}}$, then 
\[
(L_{g} u )_{\bar{j}}^{\ i}= = g^{i\bar{k}} u_{, \bar{k}\bar{j}} = \frac{\partial X^i}{\partial\bar{z}^j},
\]
hence the $L^2$ norm is 
\[
|L_{g} u|^2 = g_{\alpha\bar{\beta}}g^{\mu\bar{\lambda}}\frac{\partial X^{\alpha}}{\partial\bar{z}^{\lambda}}\overline{\frac{\partial X^{\beta}}{\partial\bar{z}^{\mu}}}
= |\frac{\partial X}{\partial\bar{z}}|^2_g.
\]
now we choose a fixed smooth background metric $h$ to estimate
\[
|\frac{\partial X}{\partial\bar{z}}|^2_h = h_{\alpha\bar{\beta}}h^{\mu\bar{\lambda}}\frac{\partial X^{\alpha}}{\partial\bar{z}^{\lambda}}\overline{\frac{\partial X^{\beta}}{\partial\bar{z}^{\mu}}}
\]
\[
=  h_{\alpha\bar{\beta}}h^{\mu\bar{\lambda}} g^{\alpha\bar{\eta}}u_{,\bar{\eta}\bar{\lambda}} g^{\gamma\bar{\beta}}u_{, \gamma\mu}
= \frac{1}{\Lambda_{\alpha}^2} |u_{,\bar{\alpha}\bar{\lambda}}|^2
\]
\[
\leqslant \Sigma(\frac{\Lambda_{\lambda}}{\Lambda_{\alpha}})     \Sigma \frac{1}{\Lambda_{\alpha}\Lambda_{\lambda}} |u_{,\bar{\alpha}\bar{\lambda}}|^2
\]
\[
\leqslant C( tr_{g}h) |\frac{\partial X}{\partial\bar{z}}|^2_g
\]
where we compute in some normal coordinate. And correspondingly, the $L^2$ norm of $X$ can be estimated by
\[
|X|^2_h = h_{\alpha\bar{\beta}}g^{\alpha\bar{\lambda}}u_{,\bar{\lambda}} \overline{g^{\beta\bar{\eta}}u_{, \bar{\eta}}}
\]
\[
= \frac{1}{\Lambda_{\alpha}^2}|u_{,\bar{\alpha}}|^2
\]
\[
\leqslant \Sigma(\frac{1}{\Lambda_{\alpha}}) \Sigma \frac{1}{\Lambda_{\alpha}} |u_{,\bar{\alpha}}|^2
\]
\[
\leqslant  (tr_g h) |\bar{\partial} u|^2_g.
\]

 Recall that $f = \phi'' - |\partial\phi'|^2_g $ is bounded from above, then we can estimate the $L^2$ norm of $\bar{\partial}v$ as 
\[
\int_X |\frac{\partial X}{\partial\bar{z}}|^2_h \det h\leqslant C\int_X |\frac{\partial X}{\partial\bar{z}}|^2_h \frac{1}{f} \det h 
\]
\[
\leqslant\frac{C}{\epsilon} \int_X |\frac{\partial X}{\partial\bar{z}}|^2_g (tr_g h) \det g
\]
\[
\leqslant \frac{C'}{\epsilon}\int_X |\frac{\partial X}{\partial\bar{z}}|^2_g e^{-\phi_g} \leqslant C''.
\]
note $X$ is a vector in $(1,0)$ direction, which means locally its coefficients are functions. Hence its full gradient is uniformly bounded in $L^2$ norm, i.e.
\[
\int_X |\nabla X_{\epsilon} |^2_h\det h < C
\]
for some constant independent of $\epsilon$. We claim it's also $L^1$ bounded. Recall from our choice of $\epsilon$, we have 
\[
\int_X f e^{-\phi_{\epsilon}} < C_1 \epsilon,
\]
then we can estimate 
\[
\int_X e^{F_{\epsilon}} \det h = \frac{1}{\epsilon} \int_X f e^{-\phi_{\epsilon}} < C_1,
\]
hence
\[
(\int_X |X|_h \det h )^2 \leqslant C ( \int_X |X|^2_h e^{F_{g}}\det g)^2
\]
\[
\leqslant  C (\int_X |X|^2_h (\det g)^2 e^{F_g})(\int_X e^{F_g}) 
\]
\[
\leqslant C' (\int_X |\bar{\partial}u |^2_g e^{-\phi_g}) < C''.
\]
Hence it's uniformly $L^1$ bounded, then by Poinc\'are inequality, we know $||X||_{L^2} < C $ for some uniform constant. These together imply the sequence of vector fields $X_{\epsilon}$  
are uniformly $W^{1,2}$ bounded. Now by compact imbedding theorem, there exists a vector field 
$X=X^{\alpha}\frac{\partial}{\partial z^{\alpha}} \in W^{1,2}$ such that $X_{\epsilon}\rightarrow X$ in strong $L^2$ norm. 

Moreover, observe that 
\[
(\int_X |\frac{\partial X}{\partial\bar{z}}|_h e^{-\phi_g} )^2= (\int_X |\frac{\partial X}{\partial\bar{z}}|_h e^{F_g}\det g)^2
\] 
\[
\leqslant ( \int_X |\frac{\partial X}{\partial\bar{z}}|^2_h (\det g)^2 e^{F_g} )(\int_X e^{F_g})
\]
\[
\leqslant (C \int_X |\frac{\partial X}{\partial\bar{z}}|^2_g e^{-\phi_g} )(\int_X e^{F_g})
\]
\[
\leqslant C' \epsilon \int_X e^{F_g} = C' \int_X f e^{-\phi_g} 
\]
\[
< C'' \epsilon \rightarrow 0
\]
from our choice of sequence $\epsilon$. Hence $\bar{\partial}X \rightarrow 0$ in weak $L^1$ sense, but this is enough to imply $\bar\partial X = 0 $ in the sense of distributions. 
Then $X$ is in fact a holomorphic $(1,0)$ vector field on the manifolds, and we can define $v_{\infty} = X\lrcorner 1$, which is a $-K_X$ valued holomorphic $(n-1,0)$ form.

On the other hand,
for the function $u_{\epsilon}$ itself, we have
\[
\int_X |\bar{\partial}u_{\epsilon}|^2_h \det h \leqslant C \int_X |\bar{\partial} u_{\epsilon}|^2_{g_{\epsilon}} e^{-\phi_{\epsilon}}
\] 
\[
= C\Sigma_{i=1}^{N_{\epsilon}} \lambda_i^{\epsilon}|a_i(\epsilon)|^2 \leqslant C',
\]
hence $u_{\epsilon}$ has a uniform $W^{1,2}$ bound, and it converges to a function $u_{\infty}\in W^{1,2}$ in strong $L^2$ norm. Then by condition (1), the $L^2$ norm of $u_{\infty}$ is
non-trivial. Moreover, we know the equation
\[
g^{\epsilon}_{\alpha\bar{\beta}}X_{\epsilon}^{\alpha} = u(\epsilon)_{,\bar{\beta}}
\]
holds for every $\epsilon$. Now $g^{\epsilon}_{\alpha\bar{\beta}}$ is uniformly bounded from above, hence converges to $g_{\alpha\bar{\beta}}$ in weak $L^{\infty}$, where $g_{\alpha\bar{\beta}}$ is the weak $\mathcal{C}^{1,\bar{1}}$ solution of the geodesic equation. And $X_{\epsilon} \rightarrow X$ in strong $L^2$, hence by the Remark after lemma 13, we see that
the equation
\[
g_{\alpha\bar{\beta}} X^{\alpha} =\partial_{\bar{\beta}} u_{\infty}
\]
holds in the sense of $L^2$ functions. In particular, they are equal almost everywhere. 

Finally, observe that $u_{\infty}\perp \ker\bar{\partial}$, since 
\[
\int_X u_{\infty} e^{-\phi} = \lim_{\epsilon\rightarrow 0} \int_X u_{\epsilon} e^{-\phi_{\epsilon}} = 0.
\]
Hence if $v_{\infty} $ is trivial, then $\bar{\partial}u = 0$, i.e. $u\in \ker\bar{\partial}$, which implies $u = 0$, a contradiction. So $v_{\infty}$ is non-trivial too.
\end{pf}

Notice that before taking the limits, the vector field $v_{\epsilon}$ also satisfies another equation, i.e.
\[
\partial^{\phi_{\epsilon}} v_{\epsilon} = \Sigma_{i=1}^{N_{\epsilon}} \lambda_{i}^{\epsilon} a_i(\epsilon) e_i(\epsilon).
\]
the LHS converges weakly to $\partial^{\phi}v_{\infty}$, since for any smooth testing $(n,0)$ form $W$,
\[
\int_X v_{\epsilon}\wedge\overline{\bar{\partial}W} e^{-\phi_{\epsilon}} \rightarrow \int_X v_{\infty}\wedge\overline{\bar{\partial}W}e^{-\phi}
\]
and the RHS converges to $u_{\infty}$ since condition (3). And the RHS 
\[
|| \Sigma_{i=1}^{N_{\epsilon}} \lambda^{\epsilon}_i a_i(\epsilon) e_i(\epsilon) - u_{\epsilon} ||^2 \leqslant K  \Sigma_{i=1}^{N_{\epsilon}} (\lambda^{\epsilon}_i -1) |a_i(\epsilon)|^2 
\]
converges to zero. We have equality
\[
\partial^{\phi}v_{\infty} = u_{\infty}
\]
holds in the weak sense. But since both sides of above equation are $L^2$ functions, the equation actually holds as $L^2$ functions.
This reminds us that $u_{\infty}$ might be the eigenfunction of the operator $\Box_{\phi}$ with eigenvalue $1$. In fact, we have
\begin{cor}
Let $u_{\epsilon}$ be a sequence of functions satisfying condition (1) - (3) in proposition 7, then there exists a function $u_{\infty}\in W^ {1,2}$ such that
\[
u_{\epsilon}\rightarrow u_{\infty}
\] 
in strong $L^2$ sense, and $u_{\infty}$ is a nontrivial eigenfunction of the operator $\Box_{\phi_g}$ with eigenvalue $1$.
\end{cor}
\begin{pf}
First notice $u_{\infty}\in dom(\Box_{\phi_g})$. This is because $\bar{\partial} u = \omega_g\wedge v_{\infty}$, hence $u\in \mathcal{W} \subset dom(\bar{\partial})$, and $\bar{\partial}u \in dom(\bar{\partial}^*_{\phi_g})$ since $v_{\infty}$ is holomorphic. Now for any smooth testing $(n,0)$ form $W$ with value in $-K_X$, we compute
\[
\int_X \bar{\partial}^*_{\phi_g}\bar{\partial}u_{\infty} \wedge\overline{W} e^{-\phi_g} = ( \bar{\partial}^*_{\phi_g}\bar{\partial}u_{\infty}, W )_g
\]
\[
= \langle \bar{\partial}u_{\infty}, \bar{\partial}W \rangle_g
\]
\[
= \langle \omega_g\wedge v_{\infty}, \bar{\partial}W \rangle_g
\]
\[
= \int_X v_{\infty}^{\alpha}\overline{\partial_{\bar{\alpha}}W} e^{-\phi_g}
\]
\[
= ( \partial^{\phi_g}v_{\infty}, W )_g
\]
\[
= \int_X u_{\infty} \wedge\overline{W} e^{-\phi_g}.
\]
hence $\Box_{\phi_g}u_{\infty} = u_{\infty}$ as $L^2$ functions.
\end{pf}

\section{the eigenspace decomposition of $\phi'$ (the easy case)}
In this section, we shall construct a sequence of functions $u_{\epsilon}$, which could satisfy the condition $(1) - (3)$ in proposition 12 from $\phi'_{\epsilon}$, then construct a holomorphic vector 
field from there. However, we need to discuss case by case this time, i.e. let
\[
\pi_{\perp}\phi'_{\epsilon} = \Sigma_{i=1}^{+\infty} a_i(\epsilon) e_i(\epsilon),
\]
then 
\[
\Box_{\phi_{\epsilon}}(\pi_{\perp}\phi'_{\epsilon}) = \Sigma_{i=1}^{+\infty} \lambda_i^{\epsilon} a_i(\epsilon) e_i(\epsilon).
\]
Note the restriction from the vanishing of $Ding$-functional gives
\begin{equation}
\Sigma_{i=1}^{+\infty} ( \lambda_i^{\epsilon}-1) | a_i(\epsilon)|^2 < C\epsilon
\end{equation}
by passing to the chosen subsequence $\epsilon_j$. And notice that 
\[
\int_X |\bar{\partial}\phi'_{\epsilon}|^2_h  \leqslant C \int_X |\bar{\partial}\phi'_{\epsilon}|^2_{g_{\epsilon}} e^{-\phi_{\epsilon}}
\]
\[
\leqslant C \int_X \phi''_{\epsilon} e^{-\phi_{\epsilon}} < C',
\]
then there exists a function $\psi \in W^{1,2}$ such that $\phi'_{\epsilon}\rightarrow \psi$ in strong $L^2$ norm. Hence we can assume 
\begin{equation}
\frac{1}{2} < \Sigma_{i=1}^{+\infty}  | a_i(\epsilon)|^2 < 2
\end{equation}
for $\epsilon$ small enough.
\begin{rmk}
In fact ,we have $|\phi_{\epsilon} |_{\mathcal{C}^1 } < C$, hence $|| \phi_{\epsilon} ||_{W^{1,p}} < C$ for any $p$ large. Then by compact imbedding theorem, we can assume 
\[
\phi_{\epsilon}\rightarrow \phi
\]
in $\mathcal{C}^{0,\alpha}$ norm. 
\end{rmk}

 In fact, we are going to prove 
\begin{thm}
There is a holomorphic vector field $v$ on the manifolds, such that 
\[
\omega_g\wedge v = \bar{\partial}\psi
\]
where $\psi$ is the $L^2$ limit of $\phi'_{\epsilon}$ and $g$ is the $\mathcal{C}^{1,\bar{1}}$ solution of geodesic equation. Moreover, $\psi$ is a eigenfunction of the operator $\Box_{\phi_g}$ with eigenvalue $1$, i.e.
\[
\Box_{\phi_g}\psi = \psi.
\]
\end{thm}

In order to prove this theorem, we shall discuss case by case. First there are two possibilities for the convergence of eigenvalue $\lambda_i^{\epsilon}$:
\\
\\
$Case\ 1$, there exist a finite integer $k$ such that the following two things hold

i) for each $1\leqslant i \leqslant k$, $\lambda_i^{\epsilon}\rightarrow 1$ as $\epsilon\rightarrow 0$;

ii) $\lambda_{k+1}^{\epsilon}$ does not converges to $1$.
\\
\\
$Case\ 2$, for each $1\leqslant i < +\infty$, $\lambda_i^{\epsilon}\rightarrow 1$ as $\epsilon\rightarrow 0$.
\\
\\
Let's discuss $Case\ 1$ first in this section. In this case, we shall define 
\[
u_{\epsilon}: = \Sigma_{i=1}^k a_i(\epsilon)e_i(\epsilon).
\]
Notice that the divergence of $\lambda_i^{\epsilon}$ implies $\lambda_i^{\epsilon} > 1+\delta$ for some small $\delta>0$, by passing to a subsequence. Then since $\lambda_i^{\epsilon}$ is 
a non-decreasing sequence in $i$, we have for all $i>k$
\[
\lambda_i^{\epsilon} > 1+\delta
\]
for the same subsequence. Now by equation (3), we see 
\[
C\epsilon > \Sigma_{i=k+1}^{+\infty} (\lambda_i^{\epsilon} - 1) |a_i(\epsilon)|^2 
\] 
\[
\geqslant \Sigma_{i=k+1}^{+\infty} \delta  |a_i(\epsilon)|^2,
\]
hence $\Sigma_{i=k+1}^{+\infty} |a_i(\epsilon)|^2 \rightarrow 0$ when $\epsilon\rightarrow 0$. This gives condition (1), i.e.
\[
\Sigma_{i=1}^k |a_i(\epsilon)|^2 > 1/4.
\]
condition (2) is satisfied because $\lambda_{k}^{\epsilon} \rightarrow 1$ by the assumption, and condition (3) is automatically satisfied by equation (3). Hence we can generate a holomorphic 
vector field $v_{\infty}$ from proposition (12).

Moreover, we could see $ ||\pi_{\perp}\phi'_{\epsilon} - u_{\epsilon} ||_{L^2} $ converges to zero in above argument, hence we actually have 
\[
\psi = u_{\infty}
\] 
after taking the limit. And hence it's the eigenfunction of $\Box_{\phi_g}$ with eigenvalue $1$, by corollary (14). Hence we proved theorem 15 in this case.

\section{the hard case}
Now we are going to deal with $Case\ 2$, i.e. we assume 
\[
\lambda_i^{\epsilon} \rightarrow 1
\]
for each $1\leqslant i < +\infty$. Here we still subdivide it into two subcases as follows:
\\
\\
$subCase\ 1$,  for any $1< k<\infty$, the partial sum $\Sigma_{i=1}^{k-1} |a_{i}(\epsilon)|^2 \rightarrow 0$, when $\epsilon\rightarrow 0$.
\\
\\
$subCase\ 2$, there exists a finite number $K$, such that $\Sigma_{i=1}^{K-1} |a_{i}(\epsilon)|^2$ does not converge to zero.
\\
\\
Before going to the subcases, we need a lemma first
\begin{lem}
Let $e_i(\epsilon)$ be the eigenfunction of the weighted Laplacian $\Box_{\phi_{\epsilon}}$ with eigenvalue $\lambda_i^{\epsilon}$, i.e.
\[
\Box_{\phi_{\epsilon}}e_i(\epsilon) = \lambda_i^{\epsilon} e_i(\epsilon).
\]
Suppose there exists an uniform constant $C$, such that $\lambda_i^{\epsilon} <1+ C\epsilon$, then $e_i(\epsilon)$ converges to a non-trivial eigenfunction $e_i$ of the operator 
$\Box_{\phi_g}$ with eigenvalue $1$. Moreover, suppose there is another $j\neq i$, such that $\lambda_j$ satisfies the same condition, then $e_i, e_j$ are mutually orthogonal to each other.
\end{lem}
\begin{pf}
we define $u_{\epsilon}= e_i(\epsilon)$, then condition $(1)$ and $(2)$ hold automatically. And condition (3) is also satisfied because 
\[
\int_X (|\bar{\partial} u_{\epsilon}|^2_{g_{\epsilon}} - (\pi_{\perp}u_{\epsilon})^2)e^{-\phi_{\epsilon}} = (\lambda_i^{\epsilon}-1) < C\epsilon,
\]
hence by proposition (12) and corollary (14), we get 
\[
e_i(\epsilon)\rightarrow e_i
\]
in strong $L^2$ sense, where $e_i\in W^{1,2}$ is a eigenfunction of $\Box_{\phi_g}$ with eigenvalue $1$. Now for $j\neq i$, we have similar convergence and eigenfunction $e_j$, but
\[
\int_X e_i \bar{e}_j e^{-\phi_{g}} = \lim_{\epsilon\rightarrow 0}\int_X e_i(\epsilon)\overline {e_j(\epsilon)} e^{-\phi_{\epsilon}} =0
\]
by the strong $L^2$ convergence of $e_i(\epsilon)$, and $L^{\infty}$ convergence of $\phi_{\epsilon}$.
\end{pf}

Now let's begin to discuss the $subCase\ 1$. For any fixed $k$, by equation (4), we can find a large integer $N_{\epsilon, k}$ such that 
\[
\Sigma_{i=1}^{N_{\epsilon,k}} |a_i(\epsilon)|^2 \geqslant 1/4
\]
by the assumption in this subcase, for $\epsilon$ small
\[
\Sigma_{i=k}^{N_{\epsilon,k}}|a_i(\epsilon)|^2 \geqslant 1/8.
\]
but then by equation (3), 
\[
\frac{1}{8}(\lambda_k^{\epsilon} -1)\leqslant \Sigma_{i=k}^{N_{\epsilon,k}} (\lambda_i^{\epsilon} -1)|a_i(\epsilon)|^2 < C\epsilon,
\]
because the sequence $\lambda_i^{\epsilon}$ is non-decreasing. Hence we proved for each $k$, 
\[
\lambda_k^{\epsilon} < 1+ 8C\epsilon 
\]
for $\epsilon$ small enough. Now by lemma 16, we get an eigenfunction $e_k$ for each $1\leqslant i<\infty$, and they are orthogonal to each other. However, this is impossible since the eigenspace with eigenvalue $1$ of an elliptic operator $\Box_{\phi_g}$ has only finite rank. Hence the $subCase\ 1$ actually never happens.

\section{the final case}

Let's discuss $subCase\ 2$. Under the assumption in this case, we can find $K_1$, a finite integer, to be the first number such that $\Sigma_{i=1}^{K_1 -1} |a_i(\epsilon)|^2$ does not converge to zero. Then by passing to a subsequence, we can assume  $\Sigma_{i=1}^{K_1 -1} |a_i(\epsilon)|^2 > \delta_1 $ for some fixed positive number $\delta_1$. Now consider the truncated sequence
\[
\Lambda_1(\phi') = \Sigma_{i=K_1}^{+\infty} a_i(\epsilon) e_i(\epsilon).
\]
suppose there exists another integer $K_2 > K_1$, such that  $\Sigma_{i=K_1}^{K_2 -1} |a_i(\epsilon)|^2$ does not converge to zero, and then we can assume 
$\Sigma_{i=K_1}^{K_2 -1} |a_i(\epsilon)|^2 > \delta_2 $. We can repeat this argument, to find $0< K_1<K_2<K_3<\cdots$, but we claim this process will terminate in finite steps.
\begin{lem}
There exists an finite integer $n$, such that 
\[
\Sigma_{i=K_n}^{+\infty} |a_i(\epsilon)|^2 \rightarrow 0.
\]
\end{lem}
\begin{pf}
Let's define a sequence of sequence of functions $u^{(j)}_{\epsilon}$ as
\[
u^{(0)}_{\epsilon}: =\Sigma_{i=1}^{K_1-1}a_i(\epsilon)e_i(\epsilon)
\]
\[
u^{(1)}_{\epsilon}: =\Sigma_{i=K_1}^{K_2-1}a_i(\epsilon)e_i(\epsilon)
\]
\[
\cdots
\]
\[
u^{(j)}_{\epsilon}: = =\Sigma_{i=K_{j}}^{K_{j+1}-1}a_i(\epsilon)e_i(\epsilon)
\]
and so on. We now claim $u_{\epsilon}^{(j)}$ satisfying all the conditions (1) - (3) in proposition (12). Condition (1) is satisfied automatically by assumption, and condition (2) is satisfied since
$\lambda_k^{\epsilon}\rightarrow 1$ for any fixed $k$. Condition (3) is satisfied too because of equation (3), i.e.
\[
\Sigma_{i=K_j}^{K_{j+1}-1}(\lambda_i^{\epsilon} - 1)|a_i(\epsilon)|^2< C\epsilon,
\]
then by proposition (12) and corollary (14), we see there exists an non-trivial $W^{1,2}$ function $u^{(j)}$ such that
\[
u_{\epsilon}^{(j)}\rightarrow u^{(j)}
\]
in strong $L^2$ norm. And $u^{(j)}$ is a eigenfunction of operator $\Box_{\phi_g}$ with eigenvalue $1$. However, notice that $u^{j}_{\epsilon}$ and $u_{\epsilon}^{(k)}$ are mutually orthogonal,
and by the same argument used in lemma 16, this implies
\[
u^{(j)}\perp u^{(k)}
\]
for all different $j$ and $k$. Now we can find finite many such $u^{(j)}$ since they are all in the eigenspace with eigenvalue $1$ of the weighted
Laplacian operator $\Box_{\phi_g}$, hence we proved the lemma.
\end{pf}

Next we are going to complete the proof of theorem 15. Now let's define
\[
u_{\epsilon}: =\Sigma_{i=1}^{K_n-1}a_i(\epsilon) e_i(\epsilon)
\]
where $K_n$ is the number appearing in lemma 17. Now people can check the three conditions in proposition 12 are satisfied, and hence there exists a $W^{1,2}$ function $u$ such that
\[
u_{\epsilon}\rightarrow u
\]
in $L^2$ sense, and $u$ is a eigenfunction with eigenvalue $1$ of operator $\Box_{\phi_g}$, and there is a holomorphic vector field $v$ such that 
\[
\omega_g\wedge v = \bar{\partial}u.
\]
Moreover, the difference of the $L^2$ norm is
\[
||\pi_{\perp}\phi'_{\epsilon} - u_{\epsilon} ||_{L^2} = \Sigma_{i=K_n}^{+\infty}|a_i(\epsilon)|^2 \rightarrow 0 
\]
by our choice of $K_n$, hence we have
\[
\psi = u.
\]
And we complete the proof.

\begin{rmk}
If there is no any non-trivial holomorphic vector field on $X$, then {\em proposition 12} directly implies $\phi'=0$ almost everywhere on $X\times I$ from above case by case discussion. Without using {\em corollary 14}, we don not need to invoke any eigenfunction of the first eigenspace of the weighted Laplacian operator in the limit. Hence we proved uniqueness in this case.
\end{rmk}
\section{Time direction}
Up to now, we construct a holomorphic vector field $v_t$ on a fiber $X\times{t}$ for almost everywhere $t\in [0,1]$. 
And this vector field can be computed as 
\[
v_t = \omega_g \lrcorner \bar{\partial}\psi
\]
where $\phi'_{\epsilon}\rightarrow \psi$ in strong $L^2$ norm at time $t$. Notice that there are more information to use for the convergence of
$\phi'_{\epsilon}$. In fact, we know $|\phi'|, |\phi_{t\bar{z}}|$ and $|\phi_{z\bar{t}}|$ are all uniformly bounded on $X\times I$, i.e.
\[
|\phi'|_{\mathcal{C}^1}< C,
\]
then we can assume $\phi'_{\epsilon}\rightarrow \phi' \in \mathcal{C}^1(X\times I)$, in $\mathcal{C}^{0,\alpha}$ norm. Hence the two limits actually agree with each other, i.e.
\[
\psi = \phi' 
\]
as $L^2$ functions on $X$. Now the holomorphic vector field can be written as 
\[
v_t = \omega_g \lrcorner \bar{\partial} \phi'.
\]
Then we can define the following subset of the unit interval
\[
S:= \{ t\in I ; \ there\ is\ a\ holomorphic\ vector\ field\ v_t\ on\ X\times\{ t\}\ satisfying\ \omega_{g}\wedge v_t =\bar{\partial}\phi'     \}
\]
we know the set $I - S$ has measure zero. Next we are going to prove a stronger result
\begin{prop}
The subset $S$ coincides with the whole unit interval, i.e. 
\[
S = I.
\]
\end{prop}
\begin{pf}
First recall that $\phi_{\epsilon}\rightarrow \phi$ in $\mathcal{C}^{0,\alpha}(X\times I)$ norm, by the uniform bound on $\mathcal{C}^1$ norm of $\phi$. Then on each fiber $X\times\{ t\}$, the convergence still holds, i.e.
\[
\phi_{\epsilon}\rightarrow \phi
\] 
in $\mathcal{C}^{0,\alpha}(X)$, and this implies 
\[
g_{\epsilon,\alpha\bar{\beta}} \rightarrow g_{\alpha\bar{\beta}}
\]
in the sense of distribution on the fiber $X\times\{t \}$. Pick up a point $\underline{t}\in I-S$, and a sequence $t_i \in S$ such that $t_i \rightarrow \underline{t}$. Observe that the space of all
holomorphic vector fields is finite dimensional, i.e. let 
\[
\Gamma(X): = H^0(TX),
\]
then $\Gamma$ is a finite dimensional vector space. Write $v_{t_i} = X_i\lrcorner 1 $, where $v_{t_i}\in \Gamma$ is the vector field satisfying the equation in the definition of $S$. 
Observe that $v_t$ is the unique solution to the following equation
\[
\partial^{\phi_t}v_t = \Box_{\phi_t}\phi' = \pi_{\perp}\phi'
\]
under the condition $H^{0,1}(X)=0$, then the standard $L^2$ estimate(Berndtsson[5]) gives us 
\[
|| v_t ||_h \leqslant  C || \pi_{\perp} \phi' ||_h 
\]
for some fixed metric $h$ and uniform constant $C$ independent of time $t$. Consider the sequence $\{X_i \}\in H^0(TX)$, the uniform bounds on the $L^2$ norm of $X_i$ shows it must converges under the fixed metric $h$, i.e. there exists a vector field $X\in \Gamma$ such that
\[
|| X - X_i ||_h^2 \rightarrow 0.
\]
Let's write $g_{\alpha\bar{\beta}} = g_{ \alpha\bar{\beta}}(\underline{t})$ and $g_{i, \alpha\bar{\beta}} = g_{\alpha\bar{\beta}}(t_i )$, then 
\[
|| X - X_i ||_g^2 \leqslant C || X -X_i ||_h^2,
\]
hence converges to zero too. Now we claim the equation
\[
\omega_g\wedge X = \bar{\partial}\phi' 
\]
holds in the sense of distribution. Put $\chi(z)$ be any smooth compact supported testing function on $X$(we can further assume $\chi$ is supported in some coordinate chart), we fix a pair of index $\alpha, \beta$, and compute 
\[
\int_X (g_{\alpha\bar{\beta}}X^{\alpha} - g_{i,\alpha\bar{\beta}}X_i^{\alpha} )\chi(z)\det h
\]
\[
= \int_X \chi (g_{\alpha\bar{\beta}} -g_{i, \alpha\bar{\beta}})X^{\alpha} \det h + \int_X\chi (X^{\alpha} - X_i^{\alpha} )g_{i,\alpha\bar{\beta}} \det h,
\]
since $g_{i, \alpha\bar{\beta}}$ is uniformly bounded, the second term in above equation converges to zero in strong $L^2$ sense. And the first term, we can decompose it into
\[
\int_X \chi (g_{\alpha\bar{\beta}} - g_{i, \alpha\bar{\beta}} ) X^{\alpha}\det h
\]
\[
= \int_X\chi( g_{\alpha\bar{\beta}} - g_{\alpha\bar{\beta}}^{\epsilon} ) X^{\alpha}\det h -\int_X\chi (g_{i, \alpha\bar{\beta}} - g_{i, \alpha\bar{\beta}}^{\epsilon}) X^{\alpha}\det h
 + \int_X \chi (g_{i, \alpha\bar{\beta}}^{\epsilon} - g_{\alpha\bar{\beta}}^{\epsilon}) X^{\alpha}\det h,
\]
the first and second terms converge to zero as $\epsilon\rightarrow 0$, and for the third term, we integration by parts
\[
\int_X \chi (g_{i, \alpha\bar{\beta}}^{\epsilon} - g_{\alpha\bar{\beta}}^{\epsilon}) X^{\alpha}\det h = \int_X\chi_{,\bar{\beta}} (\phi^{\epsilon}_{i, \alpha}- \phi^{\epsilon}_{\alpha})X^{\alpha}\det h
\]
\[
= \int_X\chi_{,\bar{\beta}} (t_i - \underline{t})\phi'_{,\alpha}(t)X^{\alpha}\det h
\]
\[
\leqslant A |t_i - \underline{t}|
\]
where $A$ is a constant independent of $\epsilon$. Hence 
\[
\int_X \chi (g_{\alpha\bar{\beta}} - g_{i, \alpha\bar{\beta}} ) X^{\alpha}\det h \rightarrow 0
\]
as $t_i \rightarrow \underline{t}$, and we proved 
\[
g_{i,\alpha\bar{\beta}}X_i^{\alpha} \rightarrow g_{i,\alpha\bar{\beta}}X_i^{\alpha} 
\]
in the sense of distributions. But we know $\phi'_i\rightarrow \phi'$ in $\mathcal{C}^{0,\alpha}$ norm, hence $\bar{\partial}\phi'_i\rightarrow \bar{\partial}\phi'$ in the sense of distribution too.
Finally, the limit equation 
\[
g_{\alpha\bar{\beta}}X^{\alpha} = \phi'_{,\bar{\beta}}
\]
holds in distribution sense on $X\times \{\underline{t} \}$. Now since both sides in above equation are $L^{\infty}$ functions, we see the equation actually holds in the sense of $L^2$ functions 
by the same argument in $Remark\ 1$.
\end{pf}

Now it makes sense to talk about the time derivative of vector fields $v_t $ in distribution sense, i.e. on the $\mathcal{C}^{1,\bar{1}}$ geodesic, we compute in the sense of distributions
\[
\phi''_{,\bar{\beta}} = ( g_{\alpha\bar{\beta}} X^{\alpha} )',
\]
and computation implies 
\[
(g^{\alpha\bar{\lambda}}\phi'_{,\alpha}\phi'_{,\bar{\lambda}})_{,\bar{\beta}} = \phi'_{\alpha\bar{\beta}}X^{\alpha} + g_{\alpha\bar{\beta}}(X^{\alpha})'.
\]
note the RHS is in fact equal to 
\[
\nabla_{\bar{\beta}}(\phi'_{,\alpha}X^{\alpha}) = \phi'_{,\alpha\bar{\beta}}X^{\alpha}+ \phi'_{,\alpha}X^{\alpha}_{,\bar{\beta}} = \phi'_{,\alpha\bar{\beta}}X^{\alpha},
\]
here Leibniz rule makes sense since $X$ is holomorphic. Hence we get
\[
g_{\alpha\bar{\beta}}(X^{\alpha})'  =  0
\]
which is equivalent to the vanishing of $\frac{\partial}{\partial t}v_t = 0$, i.e. we have an unchanged holomorphic vector field $v$ on the geodesic.
\\
\\
We finished the proof of uniqueness theorem by taking the holomorphic vector field 
\[
\mathcal{V}:= \frac{\partial}{\partial t} - V,
\]
then it's easy to check $\mathcal{L}_{\mathcal{V}} (i\partial\bar{\partial}\phi_t) = 0$ during the flow, hence the induced the automorphism $F$ preserves the metric along the geodesic.
\\
\\

\end{document}